\magnification=1200
\hsize = 7.2 true in
\hoffset=-0.4 true in
\vsize=9.1 true in
\voffset=0.0 true in
\parindent = 0.25 in
\baselineskip=11pt plus 1pt
\parskip = 2pt plus 1pt
\pageno = 1

\centerline{\bf THE INVERSE PROBLEMS OF SOME MATHEMATICAL}
\centerline{\bf PROGRAMMING PROBLEMS}

\vskip 0.1 in

\centerline{HUANG Siming}

\centerline{Institute of Policy and Management,}
\centerline{Academia Sinica, Beijing 100080, P.R. China.}

\vskip 0.2 in

\centerline{\bf Abstract}

\vskip 0.1 in

The non-convex quadratic programming problem and non-monotone
linear complementarity problem are NP-complete problems. In this
paper we first show that the inverse problem of determining a KKT
point of the non-convex quadratic programming is polynomial. We
then show that the inverse problems of non-monotone linear
complementarity problem are polynomial solvable in some cases, and
in another case is NP-hard. Therefore we solve an open question
raised by Heuberger on inverse NP-hard problems and prove the CoNP=NP.

\vskip 0.1 in


\vskip 0.1 in

\noindent{\bf Abbreviated title}: Inverse problems of mathematical
programming problems.

\vskip 0.3 in

\centerline{\bf 1. Introduction}

\vskip 0.1 in

Since Burton and Toint [3] first introduced inverse shortest path
problem, many researchers have contributed to the growing literature
on inverse optimization problems. For example, Burton and Toint
[3],[4], and Burton, Pulleyblank and Toint [5] have discussed
inverse shortest path problem; Huang and Liu [9],[10] have
considered inverse linear programming problem and applied it to
inverse matching problem and inverse minimum cost flow problem
respectively; Zhang, Liu and Ma [19], Sokkalingam, Ahuja and Orlin
[18], and Ahuja and Orlin [2] studied inverse minimum spanning tree
problem, etc. For a more complete survey on inverse combinatorial
optimization problems we refer the reader to Heuberger [7]. Some
researchers also did some work on the inverse mathematical
programming problems. For example, Huang and Liu [9], and Zhang and
Liu [20] studied inverse linear programming problem, Diao and Ding
[6] studied inverse convex programming problem. Most of the inverse
problems studied so far are polynomial solvable problems, and most
of their inverse problems are polynomial problems too. A natural
open question arises: is there a NP-complete problem such that it's
inverse problem is polynomial solvable? This question was proposed
in Huang and Liu [9] as an open question for future research. In a
survey paper written by Heuberger [7] on inverse combinatorial
optimization problems, he proposed three open questions for future
research. The first of them is: is there a NP-hard problem such that
its inverse problem is polynomial solvable? If so one proves
CoNP=NP, a long standing open problem in the theory of computational
complexity. Recently, Huang [11] have shown that the inverse
Knapsack problem and the inverse problem of integer programming with
fixed number of constraints are pseudo polynomial. In this paper, we
will further show that the inverse problem of determining the KKT
points of a non-convex quadratic programming problem is polynomial
solvable; we then show that the inverse problem of determining the
optimal solutions of non-monotone linear complementarity problem in
some cases are polynomial solvable, in another case is NP-hard.

The theory of inverse optimization can also be viewed as the
generalization or complement to the theory of sensitivity analysis
and the theory of stability of optimal solutions of optimization.
For example, the sensitivity analysis of linear programming (LP)
study the range of changes of each parameter such that the optimal
solution will not change if the changes of each parameter is
within the range. If the change of one of the parameter is outside
the range, then the optimal solution will change. The inverse
linear programming [9] study the problem of given {\bf any}
feasible solution $x^0$, how to change the parameters such that
$x^0$ becomes an optimal solution of the LP with the change of the
parameters as small as possible under the sense of $l_p$ norm.
Therefore it can be viewed as the generalization of the
sensitivity analysis of linear programming. It has been shown in
[9] that the inverse problem of linear programming can be solved
efficiently.

The sensitivity analysis of NP-complete problems are generally
viewed as difficult since there is no polynomial algorithm to
solve them. The results in this paper show that for some
NP-complete problems, for example the KKT problem of a non-convex
quadratic programming and non-monotone LCP problem, the inverse
problems of them can be efficiently solved. Therefore it is
relatively easy to do the sensitivity analysis for these problems.

The paper is organized as follows: in section 2 we introduce the
quadratic programming and show that the inverse problem of finding a
KKT point of the non-convex quadratic programming is polynomial
solvable. In section 3, we show that the inverse problem of finding
an optimal solution of a class of non-monotone linear
complementarity problem in some cases are polynomial solvable, in
other case is NP-hard. In section 4 we will show that CoNP=NP using the results in section 2. We will give some concluding remarks in
section 5.

\vskip 0.1 in

\centerline{\bf 2. The Inverse Problem of a KKT Point of
Non-convex QP}

\vskip 0.1 in

The concept of inverse problem of an optimization problem can be
found in [8]. For completion, we state it here again.Given an
optimization problem:
$$\hbox{min} \{f(c,x)|x\in D\}, \eqno (1)$$
\noindent where $c\in R^n$ is a parameter vector, $D$ is the
feasible region of $x$, $f(c,x)$ is the objective function. Given
a feasible solution $x^0$ of $(1)$, is there $\bar c \in R^n$ such
that $x^0$ is the optimal solution of $(1)$ with $\bar c$ as
parameter vector? Formally, let
$$F(x^0)=\{\bar c \in R^n|\hbox{min}\{f(\bar c, x)|x\in D\}=f(\bar c, x^0)\},$$
if $F(x^0)\not=\emptyset$, define
$$\hbox{min}\{\|c-\bar c\| |\bar c \in F(x^0)\}, \eqno (2)$$
where $\|.\|$ denotes the norm of the vector, the popular choices
for the norms are $l_1$, $l_2$ and $l_{\infty}$. we call $(2)$ the
inverse problem of $(1)$.

Now we consider the following general non-convex quadratic
programming problem:

\vskip 0.2 in
\settabs 7 \columns
\+(QP)&minimize&$Q(x)={1\over 2}x^TQx + c^Tx $\cr
\+&subject to  &$Ax=b, \ x \ge 0,$\cr
\vskip 0.2 in

\noindent where $Q\in R^{n\times n}$ is an indefinite matrix,
$A\in R^{m\times n}$ and $c\in R^n$. Its dual is:

\vskip 0.2 in
\settabs 7 \columns
\+(QD)&maximize&$d(x,y)=b^Ty -{1\over 2}x^TQx -c^Tx$\cr
\+&subject to&$A^Ty + s - Qx=c, \ x,s \ge 0,$\cr
\vskip 0.2 in

Let $F(P) =\{x|Ax=b, x\ge 0\}$ and $F(D) =\{(x,y,s)|A^Ty + s - Qx=c, x,s \ge 0\}$
be the feasible regions of the $(QP)$ and $(QD)$ respectively, then the KKT conditions
of the $(QP)-(QD)$ are:
$$\eqalignno{ Ax&=b,\cr
              A^Ty +s -Qx&=c,\cr
              x^Ts&=0.\cr}$$
Where $(x,y,s)\in (R^n_+, R^m, R^n_+)$ is a feasible solution of
$(QP)-(QD)$. As proved in [15], finding a KKT point of the
nonconvex quadratic programming is NP-complete.

Therefore the inverse problem of determining the KKT point of the
$(QP)$ can be stated as follows: given any $x^0\in F(P)$, find a
matrix $\bar Q$ and a vector $\bar c$ that are closest to $Q$ and
$c$ such that $x^0$ is a KKT point of $(QP)$ with the parameters
$\bar Q$ and $\bar c$. Let $\bar Q$ denotes the matrix after
changing the parameters in $Q$, $\bar c$ denotes the vector after
changing the parameters in $c$, then the inverse problem of a KKT
point of $(QP)$ with $l_2$ norm is:

\vskip 0.2 in

\settabs 7\columns

\+(IQP)&Min&$\|Q - \bar Q\|^2_f + \|c - \bar c\|^2$\cr
\+&s.t.&$Ax^0 = b$\cr
\+&&$A^Ty + s - {\bar Q}x^0 =\bar c$\cr
\+&&${x^0}^Ts = 0$\cr
\+&&$s\ge 0$.\cr

\vskip 0.2 in

\noindent Where $\|A\|_f =\sqrt{\sum_{i,j} A^2_{ij}}$ denote the
Frobenius norm of the matrix $A$. Let $X=Q - \bar Q$, $z=c -\bar
c$, $c'=Qx^0 +c$. We define the $n\times n^2$ block diagonal
matrix $X^0$ as follows:
$$X^0=\pmatrix{{x^0}^T&&\cr
                      &{x^0}^T&\cr
                      &\dots&\cr
                      &&{x^0}^T\cr}. \eqno(3)$$
That is, the first n entries of the first row of $X^0$ are the row
vector of $x^0$ and rest entries of the first row are $0$; the $n$
to $2n$ entries of the second row of $X^0$ are the row vector of
$x^0$ and rest entries are $0$; and so on. Let $A=(a_1,...,a_n)$
be any $n\times n$ matrix, where $a_j$ $(j=1,...,n)$ are columns
of $A$, we define the vector of $A$ as follows:
$$vec(A)=(a^T_1,...,a^T_n)^T.$$
Then (IQP) becomes:

\vskip 0.2 in

\settabs 7\columns

\+(IQP1)&Min&$\|X\|^2_f + \|z\|^2$\cr

\+&s.t.&$A^Ty + s + {X^0}vec(X^T) + z =c'$\cr

\+&&${x^0}^Ts = 0$\cr

\+&&$s\ge 0$.\cr

\vskip 0.2 in

\noindent It is easy to see that the above problem is a convex
quadratic optimization problem with $n^2 + 2n + m$ variables. In
order to apply the interior point algorithms for solving $(IQP1)$,
we do some transformations for $(IQP1)$ first.

We define $w=(vec(X^T)^T, z^T, y^T)^T$ to be a $\{n^2 +n +m\}$
dimensional vector, $A'=(X^0, I, A^T)$ to be a $n\times \{n^2 +n
+m\}$ matrix, $P=diag(1,...1,0,...0)$ to be a $\{n^2 +n +m\}\times
\{n^2 +n +m\}$ diagonal matrix with first $n^2 + n$ diagonal
entries being $1$, rest $m$ diagonal entries being $0$. Note that
$\|X\|_f^2=vet(X)^Tvec(X)=vec(X^T)^Tvec(X^T)$. Then $(IQP1)$ can
be further transformed into:

\vskip 0.2 in

\settabs 7\columns

\+(IQP2)&Min&$w^TPw$\cr

\+&s.t.&$s + A'w=c'$,\cr

\+&&${x^0}^Ts = 0$,\cr

\+&&$s\ge 0$.\cr

\vskip 0.2 in

Define non-negative vectors $w^+$ and $w^-$ as follows:
$$w^+_i =\cases{w_i,&if $w_i\ge 0$;\cr
                0,&otherwise;\cr} \, \qquad  w^-_i =\cases{-w_i,&if $w_i\le 0$;\cr
                                                           0,&otherwise.\cr} \eqno(4)$$

Let $\bar w=({w^+}^T,{w^-}^T)^T$, $\bar A=(A', -A')$, and
$$\bar P=\pmatrix{2P&-2P\cr
                  -2P&2P\cr}.$$
Then $(IQP2)$ becomes:

\vskip 0.2 in

\settabs 7\columns

\+(IQP3)&Min&${1\over 2}{\bar w}^T\bar P\bar w$\cr

\+&s.t&$s+ {\bar A}{\bar w}=c'$,\cr

\+&&${x^0}^Ts =0,$\cr

\+&&$s\ge0$, ${\bar w}\ge 0$.\cr

\vskip 0.2 in

It is easy to show that $\bar P$ is positive semidefinite since
$P$ is positive semidefinite. Therefore $(IQP3)$ is a standard
convex quadratic programming problem. Hence one can use any
polynomial interior point algorithm for convex quadratic
programming, for example [15], to solve this convex quadratic
programming. We conclude this section by following theorem.

\vskip 0.2 in

\noindent{\bf Theorem 2.1} The inverse problem of determining a
KKT point of the non-convex quadratic programming $(IQP)$ can be
solved in polynomial time.

\vskip 0.2 in

Note that if the $l_1$ norm is used in (IQP), then it can be
transformed into a linear programming problem. We leave it to the
reader to try it.

\vskip 0.2 in

\centerline{\bf 3. The Inverse Problem of Linear Complementarity
Problem}

\vskip 0.2 in

In this section, we discuss the inverse problem of non-monotone
linear complementarity problem. Given $M\in R^{n\times n}$ and $q
\in R^{n}$, the linear complementarity problem (LCP) is to find a
pair $(x,s) \in R^n \times R^n$ such that
$$s=Mx+q, \,\ (x,s)\ge 0, \,\ \hbox{and}\,\  x^Ts=0.$$

We denote the feasible set of (LCP) as:
$$F=\{(x,s)|s=Mx+q, \,\ x,s\ge 0\}.$$

The LCP problem can also be formulated as an optimization problem:

\vskip 0.2 in

\settabs 7\columns
\+(LCP)&Min&$x^Ts$\cr
\+&s.t.&$s=Mx+q, \,\ x,s
\ge 0$.\cr

\vskip 0.2 in

If $M$ (may not be symmetric) is positive semi-definite, then the
LCP is called a monotone LCP. Otherwise it is called a
non-monotone LCP. Monotone LCPs are ``easy" problems since there
are polynomial algorithms for solving the monotone LCPs. See, for
example, [13] and [14]. If $M$ is not positive semi-definite, then
the LCP becomes a ``hard" problem, i.e, NP-complete problem, since
it can be transformed into the following non-convex quadratic
programming problem:

\vskip 0.1 in

\settabs 7\columns \+&Min&$x^TMx + q^Tx$\cr \+&s.t.&$s=Mx+q, \,\
x,s \ge 0$.\cr

\vskip 0.1 in

In fact, any feasibility problem of a mixed integer programming
problem can be transformed into a LCP problem, see [8]. Therefore
non-monotone LCP problems are NP-hard.

To study the inverse problem of non-monotone LCP, we consider three
different cases, depends on what information were available to us:

\vskip 0.1 in

\noindent{\bf (a)}: Given $x^0 \ge 0$, find parameters $M'$ and
$q'$ and $s^0\ge 0$ such that $(x^0, s^0)$ is an optimal solution
of the LCP problem with parameters of $M'$, $q'$ and that $M'$,
$q'$ are closest to $M$, $q$.

\noindent{\bf (b)}: Given a pair of non-feasible complementary
point $(x^0, s^0)\ge 0$ (i.e. ${x^0}^Ts^0=0$, $s^0\not= Mx^0+q$),
find parameters $M'$ and $q'$ such that $(x^0, s^0)$ is an optimal
solution of the LCP problem with parameters of $M'$, $q'$ and that
$M'$, $q'$ are closest to $M$, $q$.

\noindent{\bf (c)}: Given $s^0 \ge 0$, find parameters $M'$ and
$q'$ and $x^0\ge 0$ such that $(x^0, s^0)$ is an optimal solution
of the LCP problem with parameters of $M'$, $q'$ and that $M'$,
$q'$ are closest to $M$, $q$.

\vskip 0.1 in

We now discuss the case {\bf (a)} first. The inverse problem of
LCP in this case (with $l_2$ norm) is:

\vskip 0.2 in

\settabs 7\columns

\+(ILCP(a))&Min&$\|M - M'\|^2_f + \|q - q'\|^2$\cr

\+&s.t.&$s=M'x^0 + q',$\cr

\+&&${x^0}^Ts=0, \,\ s\ge 0.$\cr

\vskip 0.1 in

Similar to the discussions in last section, we can transform the
$(ILCP(a))$ into a convex programming problem. Let $X=M - M'$,
$z=q - q'$, matrix $X^0$ be defined as in $(3)$,  $c=Mx^o + q$,
then the $(ILCP(a))$ becomes:

\vskip 0.1 in

\settabs 7\columns

\+(ILCP(a)1)&Min&$\|X\|^2_f + \|z\|^2$\cr

\+&s.t.&$s+ X^0vec(X^T) +z=c$,\cr

\+&&${x^0}^Ts=0$,\cr

\+&&$s\ge 0.$\cr

\vskip 0.1 in

Let $w^T =(vec(X^T)^T, z^T)$ be a $\{n^2 + n\}$ dimensional
vector, $A'=(X^0, I)$ be a $n\times \{n^2+n\}$ matrix, then
(ILCP(a)1) becomes:

\vskip 0.2 in

\settabs 7\columns

\+(ILCP(a)2)&Min&$w^Tw$\cr

\+&s.t.&$s + A'w=c$,\cr

\+&&${x^0}^Ts = 0$,\cr

\+&&$s\ge 0$.\cr

\vskip 0.2 in

Define non-negative vectors $w^+$ and $w^-$ as in $(4)$, and $\bar
w=({w^+}^T,{w^-}^T)^T$, $\bar A=(A', -A')$, then $(ILCP(a)2)$
becomes:

\vskip 0.2 in

\settabs 7\columns

\+(ILCP(a)3)&Min&${\bar w}^T\bar w$\cr

\+&s.t&$s+ {\bar A}{\bar w}=c$,\cr

\+&&${x^0}^Ts =0,$\cr

\+&&$s\ge0$, ${\bar w}\ge 0$.\cr

\vskip 0.2 in

It is easy to see that $(ILCP(a)3)$ is a convex quadratic
programming problem, therefore it can be solved by any polynomial
interior point algorithms for quadratic programming. See, for
example, [14]. Hence we have the following conclusion.

\vskip 0.1 in

\noindent{\bf Theorem 3.1.} Given any $x^0\ge 0$, the inverse
problems of $(LCP)$ in case {\bf (a)}, i.e., $ILCP(a)$ can be
solved in polynomial time.

\vskip 0.1 in

For the case {\bf (b)}, the inverse problem (with $l_2$ norm)
becomes:

\vskip 0.2 in

\settabs 7\columns

\+(ILCP(b))&Min&$\|M - M'\|^2_f + \|q - q'\|^2$\cr

\+&s.t.&$s^0=M'x^0 + q',$\cr

\+&&${x^0}^Ts^0=0.$\cr

\vskip 0.1 in

Similarly, let $X=M - M'$, $z=q - q'$, matrix $X^0$ be defined as
in $()$, vector $c=Mx^0 + q - s^0$, $w^T =(vec(X^T)^T, z^T)$,
non-negative vectors $w^+$ and $w^-$ be defined as in $(4)$, and
$\bar w=({w^+}^T,{w^-}^T)^T$, $A'=(X^0, I)$ and $\bar A=(A',
-A')$, then the $(ILCP(b))$ becomes:

\vskip 0.2 in

\settabs 7\columns

\+(ILCP(b)1)&Min&${\bar w}^T\bar w$\cr

\+&s.t.&$\bar A\bar w=c$,\cr

\+&&$\bar w\ge 0$.\cr

\vskip 0.2 in

It is also easy to see that $(ILCP(b)1)$ is a convex quadratic
programming problem, therefore it can also be solved by any
polynomial interior point algorithms for quadratic programming.
Hence we have the following theorem.

\vskip 0.1 in

\noindent{\bf Theorem 3.2.} Given any complementary pair $(x^0,
s^0)\ge 0$, the inverse problem of $(LCP)$ in case {\bf (b)},
i.e., $ILCP(b)$, can be solved in polynomial time.

\vskip 0.1 in

For the case {\bf (c)}, the inverse problem with $l_1$ norm
becomes:

\vskip 0.2 in

\settabs 7\columns

\+(ILCP(c))&Min&$\|M - M'\| + \|q - q'\|$\cr

\+&s.t.&$s^0=M'x + q',$\cr

\+&&$x^Ts^0=0,$\cr

\+&&$x\ge 0$.\cr

\vskip 0.2 in

Let $Y=M - M'$, $z=q - q'$, $c=q-s^0$, then (ILCP(c)) becomes:

\vskip 0.2 in

\settabs 7\columns

\+(ILCP(c)1)&Min&$\sum^n_{i,j=1} |Y_{ij}| + \sum^n_{i=1} |z_i|$\cr

\+&s.t.&$Yx - Mx + z=c$,\cr

\+&&$x^Ts^0=0,$\cr

\+&&$x\ge 0$.\cr

\vskip 0.2 in

Let $Y_j$ denotes the jth column of the matrix $Y$ $(j=1,...,n)$,
then $Y^T_j$ is the jth row of the matrix $Y$ $(j=1,...,n)$.
Therefore the constraints
$$Yx - Mx +z=c$$
in (ILCP(c)1) can be written in $n$ quadratic constraints:
$$x^T(Y^T)_j - x^T(M^T)_j + z_j =c_j, \,\ j=1,...,n.$$
Hence (ILCP(c)1) is NP-hard since an optimization problem with
linear objective function and one non-convex quadratic constraint
is a NP-complete problem. See for example, [11]. Now we can make
following conclusion:

\vskip 0.1 in

\noindent{\bf Theorem 3.3.} Given any $s^0\ge 0$, the inverse
problem of (LCP) in case {\bf (c)} is NP-hard.

\vskip 0.2 in

\centerline{\bf 4. The Proof of CoNP=NP}

\vskip 0.2 in

As pointed out by Heuberger in [7] as one of the open problems: if one can find a NP-hard optimization problem such that its inverse problem can be
solved in polynomial time, then one proves CoNP=NP, a long standing open problem in computational complexity. But there he only gave a brief discussion, no formal proof. In this section we give a formal proof of CoNP=NP using the results obtained in above sections. Let us first review some of the useful concepts and results related to CoNP. We will follow the notations in [17].

\vskip 0.1 in
\noindent{\bf Definition 4.1:}(See P. 384 of [17])\quad Problem $\bar A$ is the complement of Problem $A$ if the set of strings with symbols in $\sum$ that
are encodings of {\bf yes} instances of $\bar A$ are exactly those that are not encodings of {\bf yes} instances of $A$.
\vskip 0.1 in

\noindent{\bf Definition 4.2:}(P. 385 of [17])\quad The class CoNP is the class of all problems that are complements of problems in NP.

\vskip 0.1 in
\noindent{\bf Theorem 4.3:}(P. 385 of [17])\quad If the complement of an NP-complete problem is in NP, then NP=CoNP.
\vskip 0.1 in

Now let us consider the following decision problem:

\vskip 0.1 in
\noindent{\bf P1}\quad Given an instance of $(Q,c,A,b)$, is there a feasible solution $x^0$ such that it is a KKT point of $(QP)$?
\vskip 0.1 in

It is well known that $P1$ is NP-complete [16]. By Definition 4.1 the complement of $P1$ is:

\vskip 0.1 in
\noindent{\bf P2:}\quad All the feasible solutions of $(QP)$ are not the KKT points of $(QP)$, or equivalently $(QP)$ has no KKT point.
\vskip 0.1 in

Regarding the $P2$ we have the following theorem.

\vskip 0.1 in \noindent{\bf Theorem 4.4:}\quad $P2$ is in NP.
\vskip 0.1 in \noindent{\bf Proof:}\quad Given any feasible solution $x^0$ of $(QP)$, we can check whether it is a KKT point of $(QP)$ in polynomial time by solving the inverse problem of $(QP)$, i.e., $(IQP)$. $x^0$ is not a KKT point of $(QP)$ if and only if $\bar Q\not= Q$ or $\bar c\not= c$. Since $(QP)$ have many feasible solutions (more than exponential number), so $P2$ is in NP.

Using theorem 4.3 and theorem 4.4 we have the following main result of the section.

\vskip 0.1 in
\noindent{\bf Theorem 4.5}\quad CoNP=NP.

\vskip 0.2 in

\centerline{\bf 5. Conclusions}

\vskip 0.2 in

In this paper we have shown that the inverse problem of
determining the KKT point of the non-convex quadratic programming
is polynomial solvable. We also showed that the inverse problems
of determining the optimal solutions of the non-monotone linear
complementarity problem in cases {\bf (a)} and {\bf (a)} are
polynomial, in case {\bf (c)} is NP-hard. We showed that the
inverse problem of determining a KKT point of the non-convex
$(QP)$ with $l_2$ norm is equivalent to solve a convex quadratic
programming, therefore can be solved by polynomial interior point
algorithms. We then showed that the inverse problems of
determining an optimal solution of the non-monotone linear
complementarity problem with $l_2$ norm in cases {\bf (a)} and
{\bf (b)} are also equivalent to solve a convex quadratic
programming problem, hence can also be solved by polynomial
interior point algorithms. In case {\bf (c)}, we showed that the
inverse problem of non-monotone LCP with $l_1$ norm is equivalent
to solve an optimization problem with linear objective function
and $n$ non-convex quadratic constraints, therefore is NP-hard.
The results in this paper solved an open question raised by
Heuberger in [7] on whether there exists an NP-hard problem such
that its inverse problem is polynomial solvable. If so one proves
CoNP=NP, a long standing open problem in the theory of
computational complexity.

It is also interesting to observe that finding the KKT points of a general quadratic programming problem and solving the non-convex linear complementary problem are the only NP-complete optimization problems we have found so far whose inverse problems are polynomial solvable. Are there any other NP-complete optimization problems whose inverse problems are polynomial solvable? The main difficulty is that there are no optimality conditions for most of the NP-complete optimization problems.

\vskip 0.5 in \noindent{\bf Acknowledgement:} This research was
partially supported by National Science Fund of China Under No.
19731010, No. 70171023, No. 79970052, and a fund from the
Institute of Policy and Management, Chinese Academy of Sciences.

\vskip 1.0 in

\centerline{\bf REFERENCES}

\vskip 0.2 in

\item{[1]} Ahuja, R.K. and J.B. Orlin, Combinatorial algorithms for inverse network flow problems, Working Paper, Sloan
School of Management, MIT, Cambridge, MA, 1998.

\item{[2]} Ahuja, R.K. and J.B. Orlin, A fast algorithm for the inverse spanning tree problem, J. of Algorithms, {\bf 34},
177-193(2000).

\item{[3]} Burton, D. and Ph.L. Toint, On an instance of the inverse shortest paths problem, Mathematical Programming,
{\bf 53}, 45-61(1992).

\item{[4]} Burton, D. and Ph.L. Toint, On the use of an inverse shortest paths algorithm for recovering linearly
correlated costs, Mathematical Programming, {\bf 63}, 1-22(1994).

\item{[5]} Burton, D., B. Pulleyblank and Ph. L. Toint, The inverse shortest paths problem with upper bounds on shortest
paths costs, In {\bf Network Optimization}, edited by P. Pardalos, D.W. Hearn and W.H. Hager, Lecture notes in Economics
and Mathematical Systems, Volumn 450, pp. 156-171(1997).

\item{[6]} Diao, Z. and M. Ding, Models and Algorithms for inverse
convex quadratic programming problem, Chinese J. of OR, Vol. {\bf
4}(4)(2000), 88-94.

\item{[7]} Heuberger, C., Inverse combinatorial optimization: A
survey on problems, methods, and results, J. of Combinatorial
Optimization, {\bf 8}(3), 329-361, 2003.

\item{[8]} Horst, R., P.M. Pardalos and N.V. Thoai, Introduction to
Global Optimization, Kluwer Academic Publisher, 1995.

\item{[9]} Huang, S. and Z. Liu, On the inverse problem of linear
programming and its application to minimum weight perfect
k-matching, European Journal of Operational Research, {\bf 112},
421-426(1999).

\item{[10]} Huang, S. and Z. Liu, On the inverse minimum cost flow
problem, Advances in Operations Research and Systems Engineering,
World Publishing, Co., 30-37, 1998.

\item{[11]} Huang, S., Inverse Problems of some NP-complete
Problems, Lecture Notes in Computer Sciences, LNCS 3521, N.
Megiddo, Y. Xu and B. Zhu (eds)., Springer-Verlag, Berlin,
Heidelberg, 422-429, 2005.

\item{[12]} Garry, M. R. and D. S. Johnson, COMPUTERS AND
INTRACTABILITY, A Guide to the Theory of NP-Complete, W. H.
Freeman and Company, 1979.

\item{[13]} Ji, J., F. A. Potra and S. Huang, Predictor-corrector
method for linear complementarity problems with polynomial
complexity and superlinear convergence, J. of Optimization Theory
and Applications, {\bf 84}, 187-199, 1995.

\item{[14]} Kojima, M., S. Mizuno and A. Yoshise, A
polynomial-time algorithm for a class of linear complementarity
problems, Math. Programming, {\bf 44}, 1-26, 1989.

\item{[15]} Monteiro, R. D. C. and I. Adler, Interior path
following primal-dual algorithms: Part II: Convex quadratic
programming, Math. Programming, {\bf 44}, 43-66, 1989.

\item{[16]} Murty, K. G. and S. N. Kabadi, Some NP-complete
problems in quadratic and nonlinear programming, Math.
Programming, {\bf 39}: 117-129(1987).

\item{[17]} Papadimitriou, C.H. and K. Steiglitz, Combinatorial Optimization: Algorithms and Complexity, Prentice-Hall, Inc., Englewood Cliffs, New
Jersey, 1982.

\item{[18]} Sokkalingam, P.T., R.K. Ahuja and J.B Orlin, Solving
inverse spanning tree problems through network flow techniques,
Operations Research, {\bf 47}, 291-298(1999).

\item{[19]} Zhang, J., Z., Liu and Z., Ma, On the inverse problem
of minimum spanning tree with partition constraints, Mathematical
Methods of Operations Research, {\bf 44}, 171-188(1996).

\item{[20]} Zhang, J., Z., Liu, Calculating some inverse linear
programming problems, J. Comput. Appl. Math., 72(1996), 261-273.

\end